
\magnification=1200
\mathsurround=1pt
\baselineskip=15pt

\nopagenumbers
\headline={\ifodd\pageno\rightheadline \else\leftheadline\fi}
\def\rightheadline{\hfil{\eightpoint\titolo}
\hfil\tenrm\folio}
\def\leftheadline{\tenrm\folio\hfil{\eightpoint\autori}
\hfil} \topskip=25pt
\def\autori{G. Dal Maso and A. Malusa}
\def\titolo{Approximation of relaxed Dirichlet problems}
\font\twelverm=cmr12
\font\twelvei=cmmi12
\font\twelvesy=cmsy10
\font\twelvebf=cmbx12
\font\twelvett=cmtt12
\font\twelveit=cmti12
\font\twelvesl=cmsl12

\font\ninerm=cmr9
\font\ninei=cmmi9
\font\ninesy=cmsy9
\font\ninebf=cmbx9
\font\ninett=cmtt9
\font\nineit=cmti9
\font\ninesl=cmsl9

\font\eightrm=cmr8
\font\eighti=cmmi8
\font\eightsy=cmsy8
\font\eightbf=cmbx8
\font\eighttt=cmtt8
\font\eightit=cmti8
\font\eightsl=cmsl8

\font\sixrm=cmr6
\font\sixi=cmmi6
\font\sixsy=cmsy6
\font\sixbf=cmbx6

\catcode`@=11 
\newskip\ttglue

\def\twelvepoint{\def\rm{\fam0\twelverm}
\textfont0=\twelverm  \scriptfont0=\ninerm
\scriptscriptfont0=\sevenrm
\textfont1=\twelvei  \scriptfont1=\ninei  \scriptscriptfont1=\seveni
\textfont2=\twelvesy  \scriptfont2=\ninesy
\scriptscriptfont2=\sevensy
\textfont3=\tenex  \scriptfont3=\tenex  \scriptscriptfont3=\tenex
\textfont\itfam=\twelveit  \def\it{\fam\itfam\twelveit}%
\textfont\slfam=\twelvesl  \def\sl{\fam\slfam\twelvesl}%
\textfont\ttfam=\twelvett  \def\tt{\fam\ttfam\twelvett}%
\textfont\bffam=\twelvebf  \scriptfont\bffam=\ninebf
\scriptscriptfont\bffam=\sevenbf  \def\bf{\fam\bffam\twelvebf}%
\tt  \ttglue=.5em plus.25em minus.15em
\normalbaselineskip=15pt
\setbox\strutbox=\hbox{\vrule height10pt depth5pt width0pt}%
\let\sc=\tenrm  \let\big=\twelvebig  \normalbaselines\rm}

\def\tenpoint{\def\rm{\fam0\tenrm}
\textfont0=\tenrm  \scriptfont0=\sevenrm  \scriptscriptfont0=\fiverm
\textfont1=\teni  \scriptfont1=\seveni  \scriptscriptfont1=\fivei
\textfont2=\tensy  \scriptfont2=\sevensy  \scriptscriptfont2=\fivesy
\textfont3=\tenex  \scriptfont3=\tenex  \scriptscriptfont3=\tenex
\textfont\itfam=\tenit  \def\it{\fam\itfam\tenit}%
\textfont\slfam=\tensl  \def\sl{\fam\slfam\tensl}%
\textfont\ttfam=\tentt  \def\tt{\fam\ttfam\tentt}%
\textfont\bffam=\tenbf  \scriptfont\bffam=\sevenbf
\scriptscriptfont\bffam=\fivebf  \def\bf{\fam\bffam\tenbf}%
\tt  \ttglue=.5em plus.25em minus.15em
\normalbaselineskip=12pt
\setbox\strutbox=\hbox{\vrule height8.5pt depth3.5pt width0pt}%
\let\sc=\eightrm  \let\big=\tenbig  \normalbaselines\rm}

\def\ninepoint{\def\rm{\fam0\ninerm}
\textfont0=\ninerm  \scriptfont0=\sixrm  \scriptscriptfont0=\fiverm
\textfont1=\ninei  \scriptfont1=\sixi  \scriptscriptfont1=\fivei
\textfont2=\ninesy  \scriptfont2=\sixsy  \scriptscriptfont2=\fivesy
\textfont3=\tenex  \scriptfont3=\tenex  \scriptscriptfont3=\tenex
\textfont\itfam=\nineit  \def\it{\fam\itfam\nineit}%
\textfont\slfam=\ninesl  \def\sl{\fam\slfam\ninesl}%
\textfont\ttfam=\ninett  \def\tt{\fam\ttfam\ninett}%
\textfont\bffam=\ninebf  \scriptfont\bffam=\sixbf
\scriptscriptfont\bffam=\fivebf  \def\bf{\fam\bffam\ninebf}%
\tt  \ttglue=.5em plus.25em minus.15em
\normalbaselineskip=11pt
\setbox\strutbox=\hbox{\vrule height8pt depth3pt width0pt}%
\let\sc=\sevenrm  \let\big=\ninebig  \normalbaselines\rm}

\def\eightpoint{\def\rm{\fam0\eightrm}
\textfont0=\eightrm  \scriptfont0=\sixrm  \scriptscriptfont0=\fiverm
\textfont1=\eighti  \scriptfont1=\sixi  \scriptscriptfont1=\fivei
\textfont2=\eightsy  \scriptfont2=\sixsy  \scriptscriptfont2=\fivesy
\textfont3=\tenex  \scriptfont3=\tenex  \scriptscriptfont3=\tenex
\textfont\itfam=\eightit  \def\it{\fam\itfam\eightit}%
\textfont\slfam=\eightsl  \def\sl{\fam\slfam\eightsl}%
\textfont\ttfam=\eighttt  \def\tt{\fam\ttfam\eighttt}%
\textfont\bffam=\eightbf  \scriptfont\bffam=\sixbf
\scriptscriptfont\bffam=\fivebf  \def\bf{\fam\bffam\eightbf}%
\tt  \ttglue=.5em plus.25em minus.15em
\normalbaselineskip=9pt
\setbox\strutbox=\hbox{\vrule height7pt depth2pt width0pt}%
\let\sc=\sixrm  \let\big=\eightbig  \normalbaselines\rm}

\def\twelvebig#1{{\hbox{$\textfont0=\twelverm\textfont2=\twelvesy
	\left#1\vbox to10pt{}\right.\n@space$}}}
\def\tenbig#1{{\hbox{$\left#1\vbox to8.5pt{}\right.\n@space$}}}
\def\ninebig#1{{\hbox{$\textfont0=\tenrm\textfont2=\tensy
	\left#1\vbox to7.25pt{}\right.\n@space$}}}
\def\eightbig#1{{\hbox{$\textfont0=\ninerm\textfont2=\ninesy
	\left#1\vbox to6.5pt{}\right.\n@space$}}}


%
\newcount\parn \newcount\forn \newcount\teon
\newcount\vol \newcount\pag
%
%
\parn=1
\def\newpar#1{{\global\forn=0 \global\teon=0
\global\advance\parn by 1} \bigskip {\noindent\bf \the\parn
. #1} \bigskip}  %
\def\phor{{\global\advance\forn by 1}
\leqno{(\the\parn.\the\forn)}}
\long\def\theo#1{{\global\advance\teon by 1}\medskip {\noindent\bf Theorem
\the\parn.\the\teon}\quad{\sl #1}\medskip}
\long\def\prop#1{{\global\advance\teon by 1}\medskip {\noindent\bf
Proposition \the\parn.\the\teon}\quad{\sl #1}\medskip}
\long\def\lemma#1{{\global\advance\teon by 1}\medskip {\noindent\bf
Lemma \the\parn.\the\teon}\quad{\sl #1}\medskip}
\def\ehse{{\global\advance\teon by 1}\medskip {\noindent\bf Example
\the\parn.\the\teon}\quad}
\def\ohss{{\global\advance\teon by1}\medskip {\noindent\bf Remark
\the\parn.\the\teon}\quad}
\def\dhef{{\global\advance\teon by 1}\medskip {\noindent\bf
Definition \the\parn.\the\teon}\quad}
\def\proof{{\noindent\bf Proof.}\quad}
\def\sqr#1#2{{\vcenter{\vbox{\hrule height.#2pt
													\hbox{\vrule width.#2pt height#1pt \kern#1pt
													\vrule width.#2pt}
             \hrule height.#2pt}}}}
\def\finedim{\hfill\hbox{$\sqr44$ \qquad} \par}
\def\sob#1#2{W^{#1}_{#2}(\Omega)}
\def\huz{H^1_0(\Omega)}
\def\elle#1{L^{#1}(\Omega)}

\def\intom#1{\int_{\Omega}{#1 \, dx}}

\def\norma#1#2{\|#1\|_{\lower 4pt \hbox{$\scriptstyle #2$}}}
\def\inside#1{\null\,\vcenter{\openup1\jot \mathsurround=0pt
\ialign{\strut\hfil${##}$\hfil
\crcr#1\crcr}}\,}
\def\bibart#1#2#3#4#5#6#7#8{\global\vol=#5 \global\pag=#7 \par
{\item{[{#1}]}{#2}: {#3}. {\it
#4}
{\ifnum\vol=0\else{{\bf\ #5}}\fi}
{(#6)}%
{\ifnum\pag=0\else{, {#7}--{#8}}\fi}%
. \medskip}}

\def\biblib#1#2#3#4#5#6{\par
{\item{[{#1}]}{#2}: {#3}. {#4}, {#5}, {#6}. \medskip}}

\def\pb#1#2{\cases{#1 & in $\Omega\,,$ \cr #2 & on
$\partial\Omega\,,$ \cr}}

\long\def\salta#1{\relax}

\def\capa#1#2{{\rm cap}(#1,#2)}
\def\capam#1#2#3#4{{\rm cap}^{#1}_{#2}(#3,#4)}
\def\mo{{\cal M}_0(\Omega)}
\def\erre{{\bf R}}
\def\enne{{\bf N}}
\def\rn{\erre^n}
\def\fm#1#2{F_{#1}(#2)}
\def\intmu#1#2#3{\int_{#1}{#2\,d#3}}
\def\huz{H^1_0(\Omega)}
\def\elmu#1#2{L^{#1}_{#2}(\Omega)}
\def\hm{H^{-1}(\Omega)}
\def\hu{H^1(\Omega)}
\def\infi#1{\infty_{#1}}
\def\cio{C^\infty_0(\Omega)}

\def\qih{Q_h^i}
\def\bih{B_h^i}

\def\wih{w_h^i}
\def\cih{C_h^i}
\def\tih{E_h^i}

\def\mih{\nu_h^i}
\def\gih{\lambda_h^i}
\def\wh{w_h}
\def\mh{\nu_h}
\def\gh{\lambda_h}
\def\th{E_h}
\def\ch{C_h}
\def\omh{\Omega_h}
\def\uh{u_h}

\def\qv#1{Q_{#1}}
\def\bi#1{B_{#1}}
\def\hun#1{H^1(#1)}

\def\me#1{M^\rho_r\,#1}
\def\ell#1#2{L^{#1}(#2)}
\def\rest{{\vcenter{\vbox{\hbox{\vrule width.5pt height5pt
\kern5pt}\hrule height.5pt}}}}
\def\kato#1{K^+_n(#1)}

\def\la{\lambda}
\def\ga{\gamma}
\def\media#1#2#3{{1 \over #3(#1)}\,\int_{#1}{#2\,d#3}}
\def\LL{\>\hbox{\vrule width.2pt
                \vbox to7pt{\vfill
                            \hrule width7pt height.2pt}}\>}
\mathchardef\emptyset="001F

\parn=0

{\noindent \bf Introduction}

\bigskip

The notion of ``relaxed Dirichlet problem'' was introduced in [6] to
describe the asymptotic behaviour of the solutions of classical Dirichlet
problems in strongly perturbed domains. Given a bounded open subset
$\Omega$ of
$\rn$, $n\ge2$, and an elliptic operator $L$ on $\Omega$, a relaxed
Dirichlet
problem can be written in the form
$$
\cases{Lu+\mu u = f& in $\Omega$,
\cr
u=0& on $\partial\Omega$,
\cr}
\leqno(0.1)
$$
where $f\in\hm$ and $\mu$ belongs to the space $\mo$ of all
positive Borel measures on $\Omega$ which do not charge any set of
capacity zero.

The main result concerning relaxed Dirichlet problems is the following
compactness theorem (see [6], Theorem~4.14): for every sequence
$\{\omh\}$ of open subsets of $\Omega$ there exist a subsequence, still
denoted by $\{\omh\}$, and a measure $\mu\in\mo$, such that for every
$f\in\hm$ the solutions $u_h$ of the Dirichlet problems
$$
\cases{Lu_h = f& in $\omh$,
\cr
u_h=0& on $\partial\omh$,
\cr}
\leqno(0.2)
$$
extended to $0$ on $\Omega\!\setminus\!\omh$, converge in
$\ell2\Omega$ to the
unique solution $u$ of (0.1). Moreover, the following density theorem holds
(see
[6], Theorem~4.16): for every $\mu\in\mo$ there exists a sequence
$\{\omh\}$ of open subsets of $\Omega$ such that for every $f\in\hm$
the solution $u$ of (0.1) is the limit in $\ell2\Omega$ of the sequence
$\{u_h\}$ of the solutions of~(0.2). The proof of this density theorem
provides an explicit approximation only when $\mu$ is the Lebesgue
measure, while it is rather indirect in the other cases, and does not
suggest any efficient method for the construction of the sets $\omh$.

The aim of this paper is to present an explicit approximation scheme for
the
relaxed Dirichlet problems (0.1) by means of sequences of classical Dirichlet
problems of the form (0.2). We assume that $\mu\in\mo$ is a Radon
measure.
The sets $\omh$ will be obtained by removing an array of small balls from
the
set $\Omega$. The geometric construction is quite simple. For every
$h\in\enne$
we fix a partition $\{\qih\}^{}_i$ of~$\rn$ composed of cubes with side
$1/h$, and we consider the set $I_h$ of all indices
$i$ such that $\qih\subset\subset\Omega$. For every $i\in I_h$ let $\bih$
be the
ball with the same center as~$\qih$ and radius~$1/2h$, and let~$\tih$ be
another
ball with the same center such that
$$
\capam L{}{\tih}{\bih} = \mu(\qih)\,.
$$
Finally, let $\th = \bigcup_{i\in I_h}\tih$ and
$\omh=\Omega\!\setminus\!\th$.
Note that the size of the hole~$\tih$ contained in the cube~$\qih$ depends
only on
the operator~$L$ and on the value of the measure~$\mu$ on~$\qih$.

By using a very general version of the
Poincar\'e inequality proved by P.~Zamboni~[15], we shall show that, if
$\mu$ belongs to the Kato space $\kato\Omega$, i.e., the potential
generated by $\mu$ is continuous, then the method introduced by
D.~Cioranescu and F.~Murat [4] can be applied, so that for every
$f\in\hm$ the solutions $u_h$ of the Dirichlet problems (0.2) converge
in $\ell2\Omega$ to the solution $u$ of the relaxed Dirichlet problem
(0.1). To prove that the same result holds also when $\mu$ is an
arbitrary Radon measure of the class $\mo$ we use the method of
$\mu$-capacities introduced in [6] and [3].

Finally, if $\mu$ is a Radon measure and $\mu\notin\mo$, then we prove
that our
construction leads to the approximation of the solutions of the relaxed
Dirichlet problem
$$
\cases{Lu+\mu_0 u = f& in $\Omega$,
\cr
u=0& on $\partial\Omega$,
\cr}
$$
where $\mu_0$ is the greatest measure of the class $\mo$ which is less
than or
equal to~$\mu$.

\medskip

\bigskip

{\noindent \bf Acknowledgements}

\bigskip

We wish to thank Umberto Mosco for having suggested a careful
examination of
all constructive approximation methods for relaxed Dirichlet problems. We
are also
indebted to Raul Serapioni, who drew our attention to the Poincar\'e
inequality
proved in~[15].

This work is part of the Project EURHomogenization-ERB4002PL910092 of
the
Program SCIENCE of the Commission of the European Communities, and of
the
Research Project \lq\lq Irregular Variational Problems\rq\rq of the Italian
National Research Council.

\medskip

\newpar{Notation and preliminaries}

Let $\Omega$ be a bounded open subset of $\rn$, $n \geq 2$. We shall
denote
by $\hu$ and $\huz$ the usual  Sobolev spaces, by~$\hm$ the dual space
of~$\huz$, by~$\elmu  p\mu$, $1 \leq p < \infty$ the usual Lebesgue space
with respect to the measure~$\mu$; if~$\mu$ is the Lebesgue
measure, we shall use the notation~$\elle p$.

For every
subset~$E$ of~$\Omega$ the (harmonic) capacity of~$E$ with respect
to~$\Omega$ is defined by
$$
\inf \intom{|\nabla u|^2}\,,
$$
where the infimum is taken over all functions~$u \in \huz$ such
that~$u ³ 1$ a.e.\ in a neighbourhood of~$E$.
We say that a property ${{\cal P}(x)}$, depending on a point~$x \in
\Omega$, holds quasi everywhere (q.e.) in~$\Omega$ if there exists
a set~$E\subseteq \Omega$, with~$\capa E\Omega =  0$, such that~${\cal
P}$
holds in~$\Omega \!\setminus\! E$. It is well known that every~$u \in
\hu$
admits a quasi-continuous representative, which is uniquely defined up to
a set of capacity zero (see, e.g.,~[16], Theorem 3.1.4). In the sequel we
shall always identify $u$ with its quasi-continuous representative.

By a Borel measure on $\Omega$ we mean a positive, countably
additive set
function with values in~$\overline{\erre}$ defined on
the~$\sigma$-field of all Borel subsets
of~$\Omega$; by a Radon measure on~$\Omega$ we mean a Borel
measure which is
finite on every compact subset of~$\Omega$. Finally,
by~$\mo$ we denote the
set of all positive Borel measures~$\mu$ on~$\Omega$ such
that~$\mu(E) = 0$
for every Borel set~$E \subseteq \Omega$ with~$\capa E \Omega =0$.
If~$\mu$ is a Borel measure and~$E$ is a Borel subset of~$\Omega$, the
Borel
measure~$\mu\LL E$ is defined by~$(\mu\LL E)(B) = \mu(E\cap B)$
for every Borel set~$B \subseteq \Omega$. If~$\mu$, $\nu$ are Radon
measures
and~$\nu$ has a density~$f$ with respect to~$\mu$, we shall
write~$\nu=f\mu$. For every $E \subseteq \Omega$ we denote by $\infi E$
the  measure of the class~$\mo$ defined by
$$
\infi E(B)= \cases{ 0\,, & if $\capa {B \cap E}{\Omega}=0\,,$ \cr
+ \infty\,, & otherwise. \cr} \phor
$$
We shall see later that these measures are used to express the classical
Dirichlet problems~(0.2) in the form~(0.1). This will allow us to treat
problems~(0.1) and~(0.2) in a unified way.

Another class of measures we are interested in is the Kato space.

\dhef The {\it Kato space} $\kato\Omega$ is the cone of
all positive Radon measures~$\mu$ on~$\Omega$ such that
$$
\lim_{r \to 0^+} \, \sup_{x \in \Omega} \, \int_{\Omega \cap \bi
r(x)}{G_n(y-x)\,d\mu(y)} = 0\,,
$$
where~$G_n$ is the fundamental solution of the Laplace operator~$-\Delta$
in~$\rn$, and~$B_r(x)$ denotes the open ball with center~$x$ and
radius~$r$.

\medskip

For every~$\mu \in \kato\Omega$ and for every Borel set~$A\subseteq
\Omega$ we
define
$$
\norma{\mu}{\kato A} = \sup_{x \in A}\int_{A}{|y-
x|^{2-
n}\,d\mu(y)} \,, \qquad {\rm if\ }n\geq3\,,
$$
$$
\norma{\mu}{\kato A} = \sup_{x \in
A}\int_{A}{\log \left( {\rm
diam}\,(A) \over |y-x| \right)\,d\mu(y)}+\mu(A) \,,
\qquad{\rm if\ }n=2.
$$
For every~$\mu
\in \kato\Omega$ it is easy to see that~$\norma{\mu}{\kato\Omega}
<+\infty$ and~$\norma{\mu}{\kato
A}$ tends to zero as diam($A$) tends to zero. We recall that every measure
in~$\kato\Omega$ is bounded and belongs to~$\hm$. For more details
about this subject
we refer to [1], [6], [9], [14]. We shall use in the following a
Poincar\'e inequality involving Kato measures.

\lemma{ Let~$A$ be a Borel subset of a ball~$B_R=B_R(x_0)$ such
that~${\rm diam}(A) \geq q\,R$ for some~$q \in (0,1)$, and let~$\mu \in
\kato{A}$. Then there exists a positive constant~$c$, depending only
on~$q$ and on the dimension~$n$ of the space, such that $$
\intmu{A}{u^2}{\mu} \leq c\,\norma{\mu}{\kato{A}}\,\intmu{\bi
R}{|\nabla  u|^2}{x}
$$
for every $u \in H^1_0(\bi R)$.}

\proof An inequality of this kind was proved by P.~Zamboni in the
case~$n\geq 3$, $A=B_R$, and~$\mu$ absolutely continuous with respect to
the Lebesgue measure. The same arguments can be adapted, up to minor
modifications, also to the general case. The main change
in the case~$n=2$ is the use of the inequality
$$
\intmu{B_R}{{1 \over {|x-y|\,|z-y|}}}{y} \leq c_q \left(\log \left({{{\rm
diam}(A)} \over  |x-z|}\right) +1\right) \qquad \forall x,z \in A,
$$
which can be proved by direct computation. \finedim

\medskip

Finally we need a sort of dominated convergence theorem for measures
in~$\hm$.

\lemma{Let~$\{\mu_h\}$ be a sequence of positive measures belonging
to~$\hm$ that converges to 0 in the weak$^*$ topology of measures
and suppose that there exists~$\mu \in \hm$ such that~$\mu_h \leq
\mu$. Then the sequence~$\{\mu_h\}$ converges to 0 strongly in~$\hm$.}

\proof This result could be obtained easily by using the strong
compactness of the order intervals in~$\hm$. However, we give here a
self-contained elementary proof. Let us define~$\nu_h=\mu-\mu_h$.
Clearly~$\norma{\nu_h}{\hm} \leq \norma{\mu}{\hm}$ and so, up to a
subsequence,
$\{\nu_h\}$ converges to~$\mu$ weakly in~$\hm$. The previous
inequality,
together with the lower semicontinuity of the norm, implies
that~$\norma{\nu_h}{\hm}$ converges to~$\norma{\mu}{\hm}$. This
shows
that~$\{\nu_h\}$ converges to~$\mu$ strongly in~$\hm$ and concludes the
proof of
the lemma.

\finedim

\medskip

Let $L\colon \huz \to\hm$ be a linear elliptic operator in divergence form
$$
L\,u = -{\rm div}\,(A\,\nabla u)\,,
$$
where $A = A(x) = (a_{ij}(x))$ is a symmetric~$n\times n$ matrix
of bounded measurable functions satisfying, for a suitable
constant~$\alpha >0$, the ellipticity condition
$$
\alpha\,|\xi|^2 ² \sum_{i,j = 1}^n {a_{ij}(x)\,\xi_i\,\xi_j} ²
\alpha^{-1}\,|\xi|^2
$$
for a.e. $x$ in $\Omega$, and for every~$\xi \in \rn$.

A set function ${\rm cap}^L_\mu$ can be associated with every
measure~$\mu$ in
the class~$\mo$.

\dhef Let $\mu \in \mo$. For every open set~$A
\subseteq \Omega$ and for every Borel set~$E \subseteq A$ we define
the~$\mu$-capacity of~$E$ in~$A$ corresponding to the
operator~$L$ as
$$
\capam L\mu{E}{A} = \min \left\{\langle Lu,u \rangle +
\intmu{E}{u^2}{\mu}
\,\colon\, u - 1 \in H^1_0(A) \right\},
$$
 where~$\langle \cdot,\cdot \rangle$ is the duality pairing between~$\hm$
and~$\huz$.

\medskip

The $\mu$-capacity
corresponding to~$L = -\Delta$ will be denoted by~${\rm cap}_\mu$,
while the~\hbox{$\mu$-ca}\-pacity with respect to~$\mu = \infi\Omega$
will be
denoted by~${\rm cap}^L$. The latter coincides with the classical capacity
relative to the operator~$L$ according to the definition of~[13] and
[12]. If~$L = -\Delta$ and~$\mu = \infi\Omega$, then~${\rm cap}^L_\mu$
coincides with the harmonic capacity introduced at the beginning of
this section. If~$\mu=\infi F$ for some~$F\subseteq \Omega$, and~$L$ is
any elliptic operator, then~${\rm cap}^L_\mu (E,A)={\rm cap}^L(E\cap
F,A)$ for every~$E\subseteq A$.

Some of the properties of ${\rm cap}^L_\mu$ are stated in the
following proposition.

\prop{Let $\mu \in \mo$,
$A$, $B$ open
subsets of~$\Omega$ and~$E$, $F$ subsets of~$A$. Then

\item{(i)}$\capam L\mu{\emptyset}{A} = 0$;

\item{(ii)}$E \subseteq F \Longrightarrow \capam L\mu{E}{A} \leq
\capam L\mu{F}{A}$;

\item{(iii)}$\capam L\mu{E \cup F}{A} \leq
\capam L\mu{E}{A} + \capam L\mu{F}{A}$;

\item{(iv)}$A \subseteq B \Longrightarrow \capam
L\mu{E}{A} \geq
\capam L\mu{E}{B}$;

\item{(v)}$\alpha \,\capam{}\mu{E}{A} \leq \capam L\mu{E}{A} \leq
\alpha^{-1}\,\capam{}\mu{E}{A} \leq \alpha^{-1} \capam{}{}{E}{A}$;

\item{(vi)}if~$\{ E_h\}$ is an increasing sequence of subsets of~$A$
and ${E=\cup_h E_h}$, then $
\capam{L}\mu{E}{A}=\sup_h \capam{L}\mu{E_h}{A}$.

}

\proof See [6], Proposition 3.11, Theorem 3.10 and [5], Theorem 2.9.
\finedim

\medskip

Now we introduce the notion of relaxed Dirichlet problems.
\dhef Given $\mu \in \mo$ and $f \in \hm$, we say that a
function~$u$ is a solution of the {\it
relaxed Dirichlet problem}
$$
\cases{L\, u + \mu\, u = f  & in $\Omega$,\cr
u = 0  & on $\partial\Omega$,\cr} \phor
$$
if $u \in \huz \cap \elmu 2\mu$ and
$$
\langle Lu,v\rangle + \intmu{\Omega}{u\,v}{\mu} = \langle f,v \rangle
$$
for every $v \in \huz \cap \elmu 2\mu$.

\medskip

We recall that for every $f \in \hm$ there exists a unique solution~$u$ of
problem~(1.2) (see [6], Theorem 2.4).
It is easy to see that, if $E$ is a closed set, then $u$ is a
solution
of
$$
\pb{L\, u + \infi E \, u = f}{u = 0}
$$
if and only if $u = 0$ q.e. in $E \cap \Omega$ and~$u_{|_{\Omega
\!\setminus\! E}}$
is a weak solution of the the classical boundary value problem
$$
\cases{L\, u = f & in $\Omega \!\setminus\! E$, \cr
u \in H^1_0(\Omega \!\setminus\! E). \cr}
$$

\medskip

\dhef{A sequence $\{\mu_h\}$ in $\mo$ $\gamma^L$-converges
to~$\mu \in \mo$ if, for
every~$f \in \hm$, the sequence~$\{u_h\}$ of the solutions
of the
problems
$$
\cases{L\, u_h + \mu_h \, u_h = f & in $\Omega\,$, \cr u_h = 0 & on
$\partial\Omega\,$, \cr}
$$
converges strongly in $\elle 2$ to the solution $u$ of the
problem
$$
\cases{L\, u + \mu \, u = f & in $\Omega\,$, \cr u = 0 & on
$\partial\Omega\,$. \cr}
$$}

\smallskip

With every $\mu \in \mo$ we associate the lower
semicontinuous quadratic functional on~$\huz$ defined by
$$
\fm \mu u =
\langle Lu,u\rangle + \intmu{\Omega}{u^2}{\mu}\,.
$$

The following theorem shows the connection
between~$\ga^L$-convergence of the measures~$\mu_h$
and~$\Gamma$-convergence of the
corresponding functionals~$F_{\mu_h}$.

\theo{A sequence $\{\mu_h\}$ in~$\mo$
$\ga^L$-converges to the measure~$\mu \in \mo$, if and only if the
following conditions are satisfied for every~$u \in \huz$:
\item{(a)}for every sequence $\{u_h\}$ in~$\huz$ converging  to~$u$ in
$\elle2$ $$
F_\mu(u) \leq \liminf_{h \to \infty}{F_{\mu_h}(u_h)}\,;
$$

\item{(b)}there exists a sequence~$\{u_h\}$ in~$\huz$ converging
to~$u$ in~$\elle2$ such that
$$
F_\mu(u) = \lim_{h \to \infty}{F_{\mu_h}(u_h)}.
$$
}
\proof See [2], Proposition 2.9. \finedim

\medskip

Our definition of~$\ga^L$-convergence coincides with the definition
considered in~[5]. As shown in [2], Proposition~2.8, if
properties (a) and (b) hold on~$\Omega$, then they also hold for every
open set~$\Omega^\prime \subseteq \Omega$. Conversely, if (a) and (b)
hold for every open set~$\Omega^\prime \subset \subset \Omega$,
then they hold on~$\Omega$. So our definition
of~$\ga^L$-convergence
differs from the definition given in~[3] only in the fact that
now the ambient space is~$\Omega$ instead of~$\rn$.
When $L=-\Delta$, our definition coincides with the definition given
in~[6].

\ohss Let~$\{\la_h\}$ and $\{\mu_h\}$ be two sequences in~$\mo$ which
$\ga^L$-converge to~$\la$ and~$\mu$, respectively. If~$\la_h \leq
\mu_h$ for every~$h$, by Theorem 1.8 we have~$\intmu
{\Omega}{u^2}{\la} \leq \intmu
{\Omega}{u^2}{\mu}$ for every~$u \in \huz$. In particular,
if~$\mu$ is a Radon measure, then~$\la \leq \mu$.

\medskip

We briefly recall some properties of the~$\gamma^L$-convergence of
measures in~$\mo$.

\theo{For every sequence~$\{\mu_h\}$ in~$\mo$ there exists a
subsequence~$\{\mu_{h_k}\}$ which~\hbox{$\gamma^L$-con}\-verges to a
measure~$\mu$ in~$\mo$.}

\proof The proof for the case~$L = -\Delta$,
can be found in~[7], Theorem 4.14. The proof in the general case is
similar. \finedim

\medskip

\theo{Let $\{\mu_h\}$ be a sequence in~$\mo$
which~$\gamma^L$-converges to a
measure~$\mu$ in~$\mo$. Then
$$
\capam L\mu{A}{B} \leq \liminf_{h \to \infty}\,\capam
L{\mu_h}{A}{B},
$$
for every pair of open sets $A$, $B$, with~$A\subseteq B \subseteq
\Omega$.}

\proof See [5], Proposition 5.7. \finedim

\medskip

We consider now a sufficient condition for the $\gamma^L$-convergence of
a
sequence of measures of the form~$\{\infi\th\}$, where~$\{\th\}$ is a
sequence of compact subsets of~$\Omega$. In this case, if~$\omh = \Omega
\!\setminus\! \th$, the solution $\uh$ coincides with the solution of the
classical problem
$$
\cases{L\, \uh = f & in $\omh$, \cr\uh =0 & on $\partial
\Omega\,,$\cr}  $$
prolonged to zero outside~$\omh$.

Assume that $\{\th\}$ satisfies the following hypotheses, studied by
D.~Cioranescu and F.~Murat: there exist a measure~$\mu \in \sob{-
1,\infty}{}$, a
sequence $\{\wh\}$ in $\hu$, and two sequences of positive measures
of~$\hm$,
$\{\mh\}$ and~$\{\gh\}$, such that $$
\inside{
\wh \rightharpoonup 1 \qquad \hbox{weakly in } \hu\,, \cr
 \wh = 0\qquad \hbox{q.e. in }\th\,, \cr
L\,\wh = \mh - \gh\,, \cr
\mh \to \mu \qquad \hbox{strongly in } \hm\,, \cr
\gh \rightharpoonup \mu \qquad \hbox{weakly in } \hm\,, \cr}
$$
and $\langle \gh, v \rangle = 0$ for every $h \in \enne$ and for
every~$v \in \huz$, with $v = 0$ q.e. in~$\th$.

Under these hypotheses the sequence $\{\uh\}$ converges weakly
in~$\huz$ to the
weak solution~$u$ of the problem $$
\cases{L\, u + \mu\,u = f & in $\Omega$, \cr
u=0 & on $\partial \Omega$ \cr}
$$
(see [4], Th\'eor\`eme 1.2).
Later, H.~Kacimi and F.~Murat pointed out
that the hypothesis~$\mu \in \sob{-1,\infty}{}$ can be replaced by~$
\mu \in \hm$ (see [10], R\'emarque 2.4).
In conclusion, using the language introduced in Definition 1.7, the following
theorem holds.

\theo{If~$\{\th\}$ satisfies the hypotheses considered above, with~$\mu
\in \hm$,
then the sequence of measures~$\{\infi\th\}$ $\ga^L$-converges to the
measure~$\mu$.}

\medskip

\newpar{The main results}

In this section we prove that for every Radon measure~$\mu \in \mo$ the
general
approximation rule outlined in the introduction provides a sequence of
measures
of the form~$\{\infi\th\}$ which~$\gamma^L$-converges to~$\mu$
according to
Definition~1.7.

To deal with the case $\mu \in \kato\Omega$, we need the following
lemmas.

\lemma{Let $U$ and $V$ be open subsets of~$\Omega$, with~$V
\subset\subset U \subset\subset\Omega$, and let~$w$ be
the~$L$-capacitary potential of~$V$ with respect to~$U$, i.e., the unique
solution of
$$
\cases{w \in H^1_0(U)\,, & $w ³ 1$ q.e. on $V$, \cr
\langle Lw,v-w \rangle ³ 0 \,, & $\forall v \in H^1_0(U)$, $v ³ 1$
q.e. on $V$. \cr}
$$
Let us extend $w$ to~$\Omega$ by setting~$w = 0$ on~$\Omega
\!\setminus\! U$.
Then~$w \in \huz$ and~$w = 1$ q.e. on~$V$. Moreover there exist two
positive Radon measures~$\gamma$ and~$\nu$  belonging to~$\hm$ such
that~${\rm supp}\,\gamma \subseteq \partial V$, ${\rm supp}\,\nu
\subseteq \partial U$, $L\,w = \gamma - \nu$ in $\Omega$, and~$\nu
(\Omega)=\ga(\Omega) = \capam{L}{}{V}{U}$.}

We call~$\gamma$ (resp.\ $\nu$) the inner (resp.\ outer) $L$-capacitary
distribution of~$V$ with respect to~$U$.

\medskip

{\noindent\bf Proof of Lemma 2.1.}\quad It is well known (see [13],
Section 3) that there exists a positive Radon measure~$\gamma \in
H^{-1}(U)$, with~${\rm supp}\,\gamma \subseteq \partial V$, such
that~$L\,w = \gamma$ in~$\Omega$ and~$\ga(\Omega) = \capam L{}VU$.
Let us consider now the following obstacle problem
$$
\cases{z \in \huz\,, & $z ³ 0$ q.e. in $\Omega \!\setminus\! U$,\cr
\langle L\,z + \gamma,v - z \rangle \geq 0 & $\forall v \in \huz$, $v
³ 0$ q.e. in $\Omega \!\setminus\! U$.  \cr}
$$
It is well known that there exists a unique solution~$z$ of this
problem, that~$z$ is a supersolution of~$L + \gamma$, i.e.,
$L\,z + \gamma = \nu \geq 0$
for some positive measure~$\nu \in \hm$, and
that~$z \leq \zeta$ for every  supersolution~$\zeta \in \hu$ of~$L +
\gamma$
with~$\zeta \geq 0$ q.e. in~$\Omega \!\setminus\! U$ (see~[11],
Section~II.6). Since $\gamma$ is a positive measure, 0 is a
supersolution  of~$L + \gamma$. Consequently $z ² 0$ q.e. in~$\Omega$.
As~$z ³
0$ q.e. in~$\Omega \!\setminus\! U$, we conclude that~$z = 0$ q.e. in
$\Omega \!\setminus\! U$, hence~$z \in H^1_0(U)$. On the other hand~$
L\,z +
\gamma = 0$ in~$U$. As~$Lw = \ga$ on~$U$, by uniqueness we can
conclude that $z = -w$ in $U$, hence in~$\Omega$. This implies~$L\,w
= \gamma - \nu$ in~$\Omega$. As~$L\,w - \gamma = 0$ in~$U$ and
in~$\Omega \!\setminus\! \overline{U}$ we conclude that~${\rm
supp}\,\nu
\subseteq \partial U$. Since~$L\,w = \gamma - \nu$ in~$\Omega$, we
have $$
\intom{A \, \nabla w \cdot \nabla \varphi} = \intmu
{\Omega}{\varphi}{\gamma}- \intmu{\Omega}{\varphi}{\nu} \qquad
\forall
\varphi \in \huz.
$$
Let $\psi$ be a
cut-off function of class~$\cio$ such that~$\psi (x) =1$
in~$\overline{U}$. Choosing~$\varphi = {\psi (w-1)}$ as
test function we obtain
$$
\intom {A \, \nabla w \cdot \nabla w \,\psi}+ \intom {A \, \nabla
w \cdot \nabla \psi
\,(w-1)}=\intmu{\Omega}{\psi\,(w-1)}{\gamma}+\intmu{\Omega}
{\psi\,(1-w)}{\nu}
$$
and, using the fact that~$w=1$ $\gamma$-a.e. in~$\Omega$ and
$\psi\,(1-w)=1$ q.e. on~${\rm supp}\,\nu$, we obtain
$ \intom {A \, \nabla w \cdot \nabla
w}=\nu (\Omega)$. As $\ga (\Omega)= \capam L{}VU= \intom {A \, \nabla
w \cdot \nabla
w}$, we conclude that~$\nu(\Omega)=\gamma(\Omega)= \capam L{}VU$.
\finedim

\medskip

Let us fix $x^0 \in \Omega$. For every~$\rho > 0$ let~$B_\rho =
B_\rho(x^0)$
and let~$\qv\rho$ be the open cube $\{x \in \rn \colon\ -\rho < x^{}_k-
x^0_k <
\rho \hbox{ for } k=1, \ldots, n\}$. If~$0 < \rho < r$ and~$\bi r \subset
\subset
\Omega$, let~$w^\rho_r$ be the~$L$-capacitary potential of~$\bi \rho$
with
respect to~$\bi r$, and let~$\nu^\rho_r$ be the corresponding
outer~$L$-capacitary distribution.

\lemma{For every~$q \in (0,1)$ there exists a constant~$c =
c(q,\alpha,n)$, independent of the operator~$L$, such that, if~$\bi r
\subset \subset \Omega$ and~$0 < \rho ² qr$, then
$$
\media{\partial\bi r}{\varphi}{\nu^\rho_r} \leq c\,\media{\partial
\bi r}{\varphi}{\nu^{qr}_r}
$$
for every $\varphi \in H^1(\qv r)$ with $\varphi \geq 0$ q.e. in~$\qv r$.}

\proof Let us fix~$q$, $\rho$, $r$, $\varphi$ as required, and let~$u
\in \huz$ be a function whose restriction to~$\bi r$ is a solution of the
Dirichlet problem
$$
\cases{L\,u = 0 & in $\bi r\,$, \cr u - \varphi\in H^1_0(B_r)\,. \cr}
$$
We may assume that~$u=\varphi$ q.e. on the
annulus~$B_R\!\setminus\!\overline
B_r$ for some~$R>r$, so that~$u=\varphi$ q.e. on~$B_R\!\setminus\! B_r$.
By De
Giorgi's theorem, we have~$u \in C^0(\bi r)$. For every $s \in (0,r)$ we
want to
prove that $$
\media{\partial \bi s}{u}{\ga^s_r} = \media{\partial
\bi r}{\varphi}{\nu^s_r}, \phor
$$
where $\ga^s_r$ is the inner~$L$-capacitary distribution associated
with~$w^s_r$. Using the symmetry of the
operator~$L$, we get
$$
\displaylines{
0 = \intmu{\bi r}{A\,\nabla u \cdot \nabla w^s_r}x
=\intmu{\Omega}{A\,\nabla
w^s_r \cdot \nabla u}x = \cr
= \intmu{\Omega}{u}{\ga^s_r} -
\intmu{\Omega}{u}{\nu^s_r} = \intmu{\partial
\bi s}{u}{\ga^s_r} - \intmu{\partial\bi r}{\varphi}{\nu^s_r}. \cr}
$$
Since $\nu^s_r(\partial\bi r) = \capam {L}{}{\bi s}{\bi r} =
\ga^s_r(\partial\bi s)$, we obtain~(2.1).

Now we remark that, by the maximum principle, $u ³ 0$ on~$\bi r$. On
the
other hand, by Harnack's inequality,
$$
\sup_{\bi{qr}} u ² c\,\inf_{\bi{qr}} u,
$$
where the constant~$c$ depends only on~$n$, $q$, $\alpha$, (see [13],
Theorem
8.1). If we apply~(2.1) with~$s = \rho$ and~$s = qr$, we obtain
$$
\displaylines{
\media{\partial\bi r}{\varphi}{\nu^\rho_r} = \media{\partial
\bi\rho}{u}{\ga^\rho_r} ² \sup_{\bi{qr}} u ² \cr
² c\,\inf_{\bi{qr}} u ² c\,\media{\partial\bi{qr}}{u}{\ga^{qr}_r} =
c\,\media{\partial\bi r}{\varphi}{\nu^{qr}_r}, \cr}
$$
and the lemma is proved. \finedim

\medskip

For every $0 < \rho < r$, with~$\bi r \subset \subset \Omega$,
let~$M^\rho_r
\colon H^1(\qv r) \to \erre$ be the linear function defined by
$$
\me{u} = \media{\partial\bi r}{u}{\nu^\rho_r}, \phor
$$
where $\nu^\rho_r$ is the outer $L$-capacitary distribution of~$\bi\rho$
with
respect to~$\bi r$.

\lemma{For every~$q \in (0,1)$ there exists a constant~$c =
c(q,\alpha,n)$ such that, if~$\qv r \subset \subset \Omega$ and~$0 <
\rho \leq qr$, then
$$
\norma{u - \me{u}}{L^2(\qv r)} \leq c\,r\,
\norma{\nabla u}{L^2(\qv r)},
$$
for every~$ u \in H^1(\qv r)$.}

\proof Let us fix $q$, $\rho$, $r$ as required. It is not
restrictive to assume~$x^0 = 0$. Let~$Q = \qv 1$ and~$B = \bi 1$. Let us
consider the operator~$L_r$ defined by~$L_r\,u = -{\rm div}\,(A_r\,\nabla
u)$,
where $A_r(y) = A(ry)$. It is easy to check that, if~$w^\rho_r(x)$ is
the~$L$-capacitary potential of~$\bi\rho$ with respect to~$\bi
r$, then~$v^\rho_r(y) = w^\rho_r(ry)$ is the~$L_r$-capacitary potential
of~$\bi{\rho/r}$ with respect to~$B$. By Lemma~2.1 we can write~$
L_r\,v^\rho_r = \la^\rho_r - \mu^\rho_r$, with~${\rm supp}\,\la^\rho_r
\subseteq \partial\bi{\rho/r}$ and~${\rm supp}\,\mu^\rho_r \subseteq
\partial
B$. We want to prove that for every~$u \in H^1(\qv r)$ we have
$$
\intmu{\partial\bi r}{u}{\nu^\rho_r} = r^{n-2}\,\intmu{\partial
B}{u_r}{\mu^\rho_r}, \phor
$$
where $u_r(y) = u(ry)$. Let us fix~$u \in H^1(\qv r)$
and let~$\psi \in C^\infty_0(\Omega)$ be a cut-off function such that~$\psi
= 1$
on~$\partial\bi r$ and~$\psi = 0$ on~$\overline{B}_\rho$. If~$\psi_r(y) =
\psi(ry)$, then
$$
\displaylines{
\intmu{\partial\bi r}{u}{\nu^\rho_r} = \intmu{\partial\bi
r}{u\,\psi}{\nu^\rho_r} = -\intmu{\bi r}{A\,\nabla w^\rho_r \cdot \nabla
(u\,\psi)}{x} = \cr
= - r^{n-2}\,\intmu{B}{A_r\,\nabla v^\rho_r \cdot
\nabla (u_r\,\psi_r)}{y}
= r^{n-2}\,\intmu{\partial B}{u_r\,\psi_r}{\mu^\rho_r} = r^{n -
2}\,\intmu{\partial B}{u_r}{\mu^\rho_r}, \cr}
$$
which proves (2.3). Taking $u = 1$ we get~$\nu^\rho_r(\partial\bi r) =
r^{n-2}\,\mu^\rho_r(\partial B)$, so that the previous equality gives
$$
\media{\partial\bi r}{u}{\nu^\rho_r} = \media{\partial
B}{u_r}{\mu^\rho_r} \phor
$$
for every~$u \in H^1(\qv r)$.
Finally, we recall that, if~$P$ is a projection from~$\hun Q$
into~$\erre$, then the following PoincarŽ inequality holds for
every~$u$ in~$\hun Q$:
$$
\norma{u - P(u)}{\ell 2Q} \leq \beta\,\norma{P}{\left({\hun
Q}\right)^\prime}\,\norma{\nabla u}{\ell 2Q},
$$
where~$\left({\hun Q}\right)^\prime$ is the dual space of~$\hun Q$ and
the
constant~$\beta$ depends only on the dimension~$n$ of the space (see
[16], Theorem 4.2.1). Applying this result to
$$
P^\rho_r(u) = \media{\partial B}{u}{\mu^\rho_r}\,,
$$
and using (2.4), we obtain
$$
\inside{\displaystyle
\norma{u - \me{u}}{\ell 2{\qv r}}^2 = r^n\,\intmu{Q}{\left( u_r -
P^\rho_r(u_r)\right)^2}{y} \leq \cr \displaystyle
\leq \beta^2\,r^n\,\left({1 \over \mu^\rho_r(\partial
B)}\,\norma{\mu^\rho_r}{\left({\hun{Q}}\right)^\prime}\right)^2
\,\intmu{Q}{|\nabla u_r|^2}{y} = \cr \displaystyle = \beta^2\,r^2\,\left({1
\over
\mu^\rho_r(\partial
B)}\,\norma{\mu^\rho_r}{\left({\hun{Q}}\right)^\prime}\right)^2\,
\intmu{\qv r}{|\nabla u|^2}{x}\,. \cr} \phor
$$
It remains to estimate ${1 \over
\mu^\rho_r(\partial
B)}\norma{\mu^\rho_r}{\left({\hun{Q}}\right)^\prime}$. By
Lemma~2.2, applied to~$L_r$, we obtain
$$
\left| \media{\partial B}{\varphi}{\mu^\rho_r} \right| ² c\,
 \media{\partial B}{|\varphi|}{\mu^{qr}_r}
$$
for every~$\varphi \in H^1(Q)$, so that
$$
{1 \over \mu^\rho_r(\partial B)}
\norma{\mu^\rho_r}{\left({\hun{Q}}\right)^\prime} ² c\,{1 \over
\mu^{qr}_r(\partial B)}\norma{\mu^{qr}_r}{\left({\hun{Q}}\right)^\prime}.
\phor
$$
By Proposition 1.5(v) and by Lemma 2.1 we have
$$
\mu^{qr}_r(\partial B) = \capam{L_r}{}{\bi q}{B} ³ \alpha \,\capam{}{}{\bi
q}{B}. \phor $$
Let $\zeta \in C^\infty_0(\rn)$ be a cut-off function such that~$\zeta = 1$
on~$\partial B$, $\zeta = 0$ on~$\overline{B}_q$, $0 ² \zeta ² 1$ on~$B$,
and~$|\nabla\zeta|² c_q = 2/(1-q)$ on~$B$. Then, using again
Proposition~1.5(v), for every~$\varphi \in H^1(Q)$ we obtain
$$
\inside{
\displaystyle\intmu{Q}{\varphi}{\mu^{qr}_r} = \intmu{\partial
B}{\varphi\zeta}{\mu^{qr}_r} =
-\intmu{B}{A_r\,\nabla v^{qr}_r \cdot \nabla (\varphi\zeta)}{y} ² \cr
\displaystyle
² c_q\, \alpha^{-1/2}\left( \capam{L_r}{}{\bi q}{B}
\right)^{1/2}\,\norma{\varphi}{H^1(Q)} ² c_q\, \alpha^{-1}\left(
\capam{}{}{\bi {q}}{B} \right)^{1/2}\,\norma{\varphi}{H^1(Q)}. \cr}
\phor
$$
 From (2.6), (2.7), (2.8) we obtain
$$
{1 \over \mu^\rho_r(\partial B)}
\norma{\mu^\rho_r}{\left({\hun{Q}}\right)^\prime} ²Êk(q,\alpha,n),
$$
which, together with (2.5), concludes the proof of the lemma. \finedim

\medskip

For every~$r>0$ let~$\hat{Q}_r$ be the cube~$\{x \in \rn \colon\ -r \leq
x^{}_k-x^0_k < r \hbox{ for } k=1, \ldots, n\}$, so that~$Q_r$ is the interior
of~$\hat{Q}_r$.

\lemma{Let $\mu$ be a measure of~$\kato\Omega$. For every~$r > 0$,
with~$\qv r \subset\subset \Omega$, let~$\rho = \rho(r) \in (0,r)$ be the
radius
such that~$\capam{L}{}{\bi\rho}{\bi r} = \mu(\hat{Q}_r)$, and let~$M_r =
M^{\rho(r)}_r$, where~$M^{\rho(r)}_r$ is the average defined in~(2.2). Then
there
exists a function~$\omega_\mu\colon \erre_+ \to \erre_+$,
with~$\displaystyle
\lim_{r \to 0^+}{\omega_\mu(r)} = 0$, such that
$$
\norma{u - M_r\,u}{L^2_\mu(\hat{Q}_r)} ²
\omega_\mu(r)\,\norma{\nabla u}{L^2(\qv
r)} \phor
$$
for every~$u \in
H^1(\Omega)$.}

\proof First of all we prove that for every~$q \in (0,1)$ there
exists~$r_q > 0$ such that~$\rho(r) ² qr$ for~$r ² r_q$. We consider only
the
case~$n ³ 3$; the case~$n = 2$ is analogous. Since~$\mu$ is a Kato
measure,
for every~$r > 0$ we have
$$
\mu(\Omega \cap \bi r)\,r^{2-n} ² \intmu{\Omega \cap \bi
r}{|y-x^0|^{2-n}}{\mu(y)} ² \psi(r)\,,
$$
where $\psi$ is an increasing function with~$\displaystyle \lim_{r
\to 0^+}{\psi(r)} = 0$. If~$\rho = \rho(r) > qr$, then, recalling
that~$\capam{}{}{\bi{qr}}{\bi r} = c_q\,r^{n-2}$, and using
Proposition~1.5(v), we obtain
$$
\alpha c_q\,r^{n-2} \leq  \alpha\,\capa{\bi\rho}{\bi r} \leq
\capam{L}{}{\bi\rho}{\bi r}=\mu(\hat{Q}_r)\,.
$$
So we can write~$\alpha c_q\,r^{n-2} ² \mu(\hat{Q}_r) ² \mu(\Omega
\cap\bi{n\,r}) ² \beta_n\,\psi(n\,r)\,r^{n-2}$.
Choosing $r_q$ such that~$\psi(n\,r_q) < \alpha c_q / \beta_n$, we
obtain a contradiction for~$r ² r_q$. Therefore, there exists~$r_q >
0$, with~$Q_{r_q} \subset \subset \Omega$, such that~$\rho(r) ² qr$ for
every~$r
² r_q$. Since~$c_q \to +\infty$ as~$q \to 1$, we can choose~$r_q$ so that
for
every~$r > 0$, with~$Q_r \subset \subset \Omega$, there exists~$q \in
(0,1)$,
with~$r ² r_q$.

Let us fix~$q \in (0,1)$. It is clearly enough to prove~(2.9) for every~$r ²
r_q$. As~$\mu \in \kato\Omega$, by Lemma 1.2, there exists a
constant~$c_n > 0$
such that, if~$\qv r \subset \subset \Omega$, then
$$
\intmu{\hat{Q}_r}{u^2}{\mu} \leq
c_n\,\norma{\mu}{\kato{\hat{Q}_r}}\,\intmu{\bi
{nr}}{|\nabla  u|^2}{x} \phor
$$
for every $u \in H^1_0(\bi{nr})$.

Let us fix a bounded extension operator~$\Pi \colon H^1(Q_1)
\to H^1_0(B_n)$, and for every~$r > 0$ let us define the extension
operator~$\Pi_r \colon H^1(Q_r) \to H^1_0(B_{nr})$ by~$(\Pi_r u)(x)=(\Pi
u_r)(x/r)$, where~$u_r(y)=u(ry)$. It is easily seen that the boundedness
of~$\Pi$ implies the existence of a constant~$k_n > 0$ such that
$$
\int_{\bi{nr}}{|\nabla(\Pi_r\,v)|^2\,dx} \leq k_n\,\left(
\int_{\qv r}{|\nabla v|^2\,dx} + {1 \over r^2}\,\int_{\qv r}{v^2\,dx}
\right) \phor
$$
for every $v \in \hun{\qv r}$. Note that, if~$v \in\hu$ and $Q_r \subset
\subset
\Omega$, then~$v=\Pi_r\,v$ q.e. on~$\hat{Q}_r$, since both functions are
quasi
continuous and coincide on~$Q_r$. Using~(2.10) and~(2.11), for every~$u
\in
H^1(\Omega)$ we obtain  $$ \displaylines{
\intmu{\hat{Q}_r}{\left(u - M_r\,u\right)^2}{\mu}
\leq c_n\,\norma{\mu}{\kato{\hat{Q}_r}}\,\intmu{\bi{nr}}
{\bigl(\nabla\left(\Pi_r(u - M_r\,u)\right)\bigr)^2}{x} \cr
\leq c_n\,k_n\,\norma{\mu}{\kato{\hat{Q}_r}}\,\left( \int_{\qv
r}{|\nabla
u|^2\,dx} + {1 \over r^2}\,\int_{\qv r}{\left(u - M_r\,u\right)^2\,dx}
\right). \cr}
$$
As $r ² r_q$, we have~$\rho = \rho(r) ² qr$, so that Lemma~2.3 implies
that
$$
{1 \over r^2}\,\intmu{\qv r}{\left(u - M_r\,u\right)^2}{x} ²
c^2\,\intmu{\qv
r}{|\nabla u|^2}{x},
$$
hence
$$
\intmu{\hat{Q}_r}{\left(u - M_r\,u\right)^2}{\mu} \leq
c_n\,k_n\,(1 + c^2)\,\norma{\mu}{\kato{\hat{Q}_r}}\,\int_{\qv r}{|\nabla
u|^2\,dx},
$$
for every $r ² r_q$ and for every~$u \in H^1(\Omega)$.
Since~$\norma{\mu}{\kato{\hat{Q}_r}}$ tends to zero as~$r$ tends to zero,
the
statement is proved. \finedim

\medskip

We are now in a position to prove our result for Kato measures.
Let~$\{\qih\}^{}_{i\in {\bf Z}^n}$ be the partition of~$\rn$ composed of the
cubes
$$
Q^i_h=\{x \in \rn \colon \ i_k/h \leq x_k < (i_k+1)/h \hbox{ for } k=1,
\ldots,n
\}.
$$

\theo{Let $\mu\in \kato\Omega$.
Let~$I_h$ be the set of all indices~$i$ such that~$\qih \subset\subset
\Omega$.
For every $i\in I_h$ let $\bih$ be the ball with the same center as~$\qih$
and
radius~$1/2h$, and let~$\tih$ be another ball with the same center such
that
$$
\capam L{}{\tih}{\bih} = \mu(\qih)\,.
$$\nobreak
Define $\th = \bigcup_{i\in I_h}\tih$. Then the measures~$\infi{\th}$
$\gamma^L$-converge to~$\mu$ as~$h \to \infty$.}

\goodbreak
\proof Let $v^i_h$ be the $L$-capacitary potential
of~$\tih$ with respect to~$\bih$, extended to~0 on~$\Omega$, and
let~$\wih = 1 -
v^i_h$. By Lemma 2.1, we obtain~$L\,\wih = \mih - \gih$ in~$\Omega$,
with~$\mih$, $\gih \in \hm$, $\mih \geq 0$, $\gih \geq 0$, ${\rm
supp}\,\mih \subseteq \partial\bih$, ${\rm supp}\,\gih \subseteq
\partial\tih$, and
$$
\mih(\qih) = \gih(\qih) = \capam L{}{\tih}{\bih} = \mu(\qih). \phor
$$
Let us define~$\wh \in \hu$ as
$$
\wh = \cases{\wih & in $\bih \!\setminus\! \tih$, \cr
0 & in $\tih$,\cr
1 & elsewhere \cr} \phor
$$
and the measures $\mh$ and $\gh$ as
$$
\mh = \sum_{i \in I_h}{\mih}, \qquad \gh = \sum_{i \in I_h}{\gih}.
\phor
$$
We want to prove that all hypotheses of Theorem~1.12 hold for~$\wh$
and~$\mh$.

\medskip

First of all, we prove that~$\wh$ converges weakly to~1 in~$\hu$. Since, by
the maximum principle, $0 \leq \wh \leq  1$
in~$\Omega$, we have that~$\{\wh\}$ is bounded
in~$\elle 2$. On the other hand,
$$
\alpha\intom{|\nabla\wh|^2} \leq \sum_{i \in I_h}{\capam
L{}\tih\bih} = \sum_{i \in I_h}{\mu(\qih)} \leq \mu(\Omega)\,.
$$
Thus $\{\wh\}$ is bounded in $\hu$ so that there exist a
subsequence (still denoted~$\{\wh\}$) and a function~$w \in \hu$, such
that~$\{\wh\}$ converges to~$w$ weakly in~$\hu$, and hence strongly
in~$\elle2$.
We are going to show that $w = 1$ a.e. in~$\Omega$, using the arguments
of
D.~Cioranescu and F.~Murat (see [4], Th\'eor\`eme 2.2). Let us consider the
family~$\{\cih\}^{}_{i \in{\bf Z}^n}$ of all open balls with radius~$(\sqrt
n - 1)/2h$ and centers in the vertices~$i/h$ of the cubes~$\qih$. In these
balls we have $\wh = 1$. Therefore, if we define $\ch$ as the union of the
balls~$\cih$ contained in~$\Omega$, we have~$\wh\,\chi_{{}_{\ch}} =
\chi_{{}_{\ch}}$,
where $\chi_{{}_{\ch}}$ is the characteristic function of~$\ch$.
Since~$\{\chi_{{}_{\ch}}\}$ converges to a positive constant in the weak$^*$
topology
of~$\elle\infty$, passing to the limit in the equality~$\wh\, \chi_{{}_{\ch}}
=
\chi_{{}_{\ch}}$ we obtain~$w = 1$ a.e.\ in~$\Omega$.

It remains to prove that the
measures~$\mh$ defined in~(2.14) converge to~$\mu$ in the
strong topology of~$\hm$. Indeed, since~$\wh$ converges to~1
weakly in~$\hu$,
this implies also that~$\gh$ converges weakly to~$\mu$ in~$\hm$.

For every $h \in \enne$ we introduce the
polyrectangle~$P_h = \bigcup_{i \in I_h}{Q}^i_h$
and we define~$S_h = \Omega \!\setminus P_h\!$. Moreover, for every
$\varphi \in
\huz$, we consider the function
$$
\varphi_h = \sum_{i \in
I_h}{\left( M^i_h\,\varphi \right)\,\chi_{{}_{\qih}}},
$$
where, according to (2.2),
$$
M^i_h\,\varphi = \media{\partial\bih}{\varphi}{\mih}\,,
$$
and we define~$\varepsilon_h =\norma{\mu
\LL S_h}{\hm}$. Note that~$\{\varepsilon_h\}$ tends to zero by
Lemma~1.3.
Re\-cal\-ling that~$\mu(\qih) = \mih(\partial\bih)$ and using the PoincarŽ
inequality~(2.9), we have that,
$$
\displaylines{
|\langle \mh,\varphi \rangle - \langle \mu,\varphi \rangle| =
\left| \sum_{i \in I_h} {\mu(\qih) \over \mih(\partial\bih)}
\int_{\partial\bih}{\varphi\,d\mih} - \sum_{i \in I_h}
\intmu{\qih}{\varphi}{\mu} - \intmu{S_h}{\varphi}{\mu} \right| \leq \cr
\leq \intmu{P_h}{|\varphi - \varphi_h|}{\mu} +
\intmu{S_h}{|\varphi|}{\mu} \leq
\left(\mu(\Omega) \intmu{P_h}{\left(\varphi - \varphi_h\right)^2}{\mu}
\right)^{1
/ 2} + \norma{\mu \LL S_h}{\hm}\,\norma{\varphi}{\huz} = \cr
= \left(\mu(\Omega)
\sum_{i \in I_h} \norma{\varphi -  M_h^i\,\varphi}{L^2_\mu(\qih)}^2
\right)^{1
/ 2}  + \varepsilon_h\,\norma{\varphi}{\huz} \leq \cr
\leq \left(\mu(\Omega) \sum_{iÊ\in I_h}
\omega({1 / h})^2\,\norma{\nabla \varphi}{L^2(\qih)}^2 \right)^{1 /
2}  + \varepsilon_h\,\norma{\varphi}{\huz} \leq
 \left(\omega({1 /
h})\,\mu(\Omega)^{1 / 2} + \varepsilon_h\right)\norma{\varphi}{\huz}.
\cr}
$$
Thus we obtain
$$
\norma{\mh - \mu}{\hm} \leq \mu(\Omega)^{1/2}\omega(1 / h) +
\varepsilon_h,
$$
hence~$\{\mh\}$ converges to~$\mu$ strongly in~$\hm$.
Therefore~$\{\infi\th\}$ $\ga^L$-converges to~$\mu$ by Theorem~1.12.
\finedim

\medskip

In order to generalize this result to every Radon measure we need
the following results.

\prop{For every Radon measure~$\mu \in \mo$ there exist a
measure~$\nu \in \kato\Omega$ and
a positive Borel function~$g \colon \Omega \to [0,+\infty]$ such
that~$\mu = g\,\nu$.}

\proof See [2], Proposition 2.5. \finedim

\medskip

\prop{Let $\la \in \mo$, let~$\mu$ be a Radon measure
in~$\mo$; for every~$x \in \Omega$ let
$$
f(x)\, =\, \liminf_{r \to 0}\,{\capam L{\la}{\bi{r}(x)}{\bi{2r}(x)} \over
\mu(\bi{r}(x))}\,.
$$
Assume that $f$ is bounded. Then~$\la$ is a Radon measure and we
have~$\la =
f\,\mu$.}

\proof See [3], Theorem 2.3. \finedim

\medskip
\eject

\prop{Let $\mu$ be a positive Radon measure on~$\Omega$. Then there
exists
a  unique
pair~$(\mu_0,\mu_1)$ of Radon measures on~$\Omega$ such that:

\item{(i)}$\mu = \mu_0 + \mu_1$;

\item{(ii)}$\mu_0 \in \mo$;

\item{(iii)}$\mu_1 = \mu \LL N$, for some Borel set~$N$
with~$\capa N \Omega =0$.

}

\proof See [8], Lemma 2.1. \finedim

\medskip

We are now in a position to prove our main result in its most general form.

\theo{Let $\mu$
be a positive Radon measure on~$\Omega$.
Let~$\{\qih\}$ and~$\{\th\}$ be defined as in Theorem~2.5. Then the
following  results hold:

\item{(i)}if $\mu$ belongs to~$\mo$, then~$\{\infi{\th}\}$
$\gamma^L$-converges to~$\mu$;

\item{(ii)}if $\mu = \mu_0 + \mu_1$, with~$\mu_0$ and~$\mu_1$ as in
Proposition~2.8, then~$\{\infi{\th}\}$ $\gamma^L$-converges to~$\mu_0$.

}

\proof If $\mu$ is a Radon measure in~$\mo$,
then, by
Proposition~2.6, $\mu = g\,\nu$,
where~$\nu \in \kato\Omega$ and~$g$ is a positive Borel function.
By Theorem~1.8, there exists a
subsequence, still denoted by~$\{\th\}$, and a measure~$\la \in
\mo$, such
that~$\{\infi{\th}\}$ $\gamma^L$-converges to~$\la$.
Let~$x \in \Omega$ and let~$r > 0$
such
that~$\bi{2r}(x) \subseteq \Omega$. We want to prove that for every Borel
set~$E\subseteq B_{2r}$
$$
\capam L{\la}{E}{\bi{2r}(x)} \leq \mu(E). \phor
$$
If~$A$ and~$A^\prime$ are two open sets such that~$A^\prime
\subset\subset A
\subseteq B_{2r}(x)$ and~$h$ is small enough we have
$$
\bigcup_{\tih \cap A^\prime \neq \emptyset} \qih
\subseteq A\,,
$$
hence, by Theorem 1.5,
$$
\displaylines{
\capam L{}{\th \cap A^\prime}{\bi{2r}(x)} \leq
\sum_{\tih \cap A^\prime \neq \emptyset}
\capam L{}{\tih}{\bi{2r}(x)} \leq \cr
\leq \sum_{\tih \cap A^\prime \neq \emptyset}
\capam L{}{\tih}{\bih}
= \sum_{\tih \cap A^\prime\neq \emptyset} \mu(\qih) \leq \mu(A). \cr}
$$
Using Theorem~1.11 we obtain,
$$
\capam L{\la}{A^\prime}{\bi{2r}(x)} \,\leq \,\liminf_{h \to \infty}\,\capam
L{}{\th \cap A^\prime}{\bi{2r}(x)} \leq \mu(A)
$$
and, as $A^\prime \nearrow A$, we obtain~$\capam{L}{\la}{A}{\bi{2r}(x)}
\leq
\mu(A)$ for every open set~$A \subseteq \bi{2r}(x)$ (see Theorem
1.5(vi)). Since~$\mu$
is a Radon measure, this inequality can be easily extended to all Borel
subsets
of~$\bi{2r}(x)$. So~(2.15) is proved. Choosing~$E=B_r(x)$ in (2.15) and
applying
Proposition~2.7, we obtain
that~$\la$ is a Radon measure and that~$\la \leq \mu$.

Define, for~$k \in \enne$, the measures~$\mu^k = g^k\,\nu$, where~$g^k(x)
=
\min (g(x),k)$. As~$\mu^k \in \kato\Omega$, by Theorem~2.5 for
every~$k$ there
exists a sequence~$\{E_{k,h}\}_h$ such that~$\{\infi{E_{k,h}}\}_h$
$\ga^L$-converges to~$\mu^k$.
Since $\mu^k \leq \mu$, the construction of
Theorem~2.5 implies that~$E_{k,h} \subseteq E_h$ for every~$h$ and~$k$.
By Remark~1.9 this implies~$\la \geq \mu^k$ for every~$k$, hence~$\la
\geq \mu$.  As
the opposite inequality has already been proved, we obtain~$\la = \mu$.
Since
the~$\ga^L$-limit does not depend on the subsequence, the whole
sequence~$\{\infi{\th}\}$ $\ga^L$-converges to~$\mu$.

Let now~$\mu$ be any Radon measure on~$\Omega$. By Proposition~2.8,
we can
write~$\mu = \mu_0 + \mu_1$, with~$\mu_0 \in \mo$ and~$\mu_1 = \mu
\LL N$,
where~$N$ is a Borel set with~$\capa N \Omega=0$. Arguing as before,
let~$\la$ be
the~$\ga^L$-limit of a  subsequence
of~$\{\infi{\th}\}$.
If~$x \in
\Omega$ and~$r > 0$ is such that~$\bi{2r}(x) \subseteq \Omega$, we have
$$
\capam L{\la}{\bi{r}(x)}{\bi{2r}(x)} = \capam L{\la}{\bi{r}(x)
\!\setminus\!
N}{\bi{2r}(x)},
$$
since $\capam {}{}{N}{\bi{2r}(x)} = 0$ (see Proposition~1.5).
Therefore~(2.15), applied with~$E = \bi{r}(x) \!\setminus\! N$, gives
$$
\capam L{\la}{\bi{r}(x)}{\bi{2r}(x)} \leq \mu(\bi{r}(x) \!\setminus\! N) =
\mu_0(\bi{r}(x)).
$$
By applying again Proposition 2.7 we obtain~$\la \leq \mu_0$.

Since~$\mu_0$ is a Radon measure of~$\mo$,
by the first part of this theorem we can construct the holes~$E_{0,h}$
such that~$\{\infi{E_{0,h}}\}$ $\ga^L$-converges to~$\mu_0$.
Since~$\mu(\qih) \geq
\mu_0(\qih)$, we have~$E_{0,h} \subseteq E_h$, hence, by Remark
1.9, $\la \geq \mu_0$. As the opposite inequality has already been
proved, we obtain~$\la = \mu_0$. As before, this implies that the whole
sequence~$\{\infi{E_h}\}$ $\ga^L$-converges to~$\mu_0$. \finedim

\vfill

\eject

{\noindent \bf References}
\bigskip
\medskip
{\ninepoint
\bibart{1}{AIZENMAN M., SIMON B.}{Brownian motion and Harnack
inequality for Schr\"o\-din\-ger operators}{Comm. Pure Appl.
Math.}{35}{1982}{209}{273}

\bibart{2}{BAXTER J.R., DAL MASO G., MOSCO U.}{Stopping times and
$\Gamma$-convergence}{Trans. Amer. Math. Soc.}{303}{1987}{1}{38}

\bibart{3}{BUTTAZZO G., DAL MASO G., MOSCO U.}{A derivation
theorem for
capacities with respect to a Radon measure}{J. Funct.
Anal.}{71}{1987}{263}{278}

\bibart{4}{CIORANESCU D., MURAT F.}{Un terme \'etrange venu
d'ailleurs I \& II}{Nonlinear
Partial Differential Equations and their applications. Coll\`ege de
France Seminar. Vol II {\rm \& } III.\ {\rm Research Notes in Mathematics,
Pitman,  London}}{60}{1982}{98}{138 \& {\bf
70} (1983),  {154--178}}

\bibart{5}{DAL MASO G.}{$\Gamma$-convergence and
$\mu$-capacities}{Ann.
Scuola Norm. Sup. Pisa Cl. Sci. (4)}{14}{1987}{423}{464}

\bibart{6}{DAL MASO G., MOSCO U.}{Wiener criteria and energy
decay for
relaxed Dirichlet problems}{Arch. Rational Mech. Anal.}{95}{
1986}{345}{387}

\bibart{7}{DAL MASO G., MOSCO U.}{Wiener criterion and
$\Gamma$-convergence}{Appl. Math. Optim.}{15}{1987}{15}{63}

\bibart{8}{FUKUSHIMA M., SATO K., TANIGUCHI S.}{On the closable part
of
pre-Dirichlet forms and the fine supports of underlying
measures}{Osaka J.
Math.}{28}{1991}{517}{535}

\bibart{9}{KATO T.}{Schr\"odinger operators with singular
potentials}{Israel J. Math.}{22}{1972}{139}{158}

\bibart{10}{KACIMI H., MURAT F.}{Estimation de l'erreur dans des
probl\`emes de
Dirichlet o\`u appara\^\i t un terme \'etrange}{Partial
Differential
Equations and the Calculus of Variations, Essays in Honor of Ennio de
Giorgi{\rm, Birkha\"user, Boston}}0{1989}{661}{696}

\biblib{11}{KINDERLEHRER D., STAMPACCHIA G.}{An introduction to
variational inequalities and their applications}{Academic Press}{New
York}{1980}

\bibart{12}{LITTMAN W., STAMPACCHIA G., WEINBERGER H.F.}{Regular
points for elliptic equations with discontinuous coefficients}{Ann.
Scuola Norm. Sup. Pisa}{17}{1963}{45}{79}

\bibart{13}{STAMPACCHIA G.}{Le probl\'eme de Dirichlet pour les
\'equations elliptiques du second ordre \`a coefficients
discontinus}{Ann. Inst. Fourier (Grenoble)}{15}{1965}{189}{258}

\biblib{14}{SCHECHTER M.}{Spectra of partial differential
operators}{North-Holland}{Amsterdam}{1986}

\bibart{15}{ZAMBONI P.}{Some function spaces and elliptic partial
differential
equations}{Matematiche}{42}{1987}{171}{178}

\biblib{16}{ZIEMER W.P.}{Weakly differentiable
functions}{Springer-Verlag}{Berlin}{1989}

}

\medskip

\bigskip

\centerline{Gianni DAL MASO and Annalisa MALUSA}

\centerline{SISSA, Via Beirut 2/4}

\centerline{34014 TRIESTE, Italy.}
\bye